\documentclass[british,authoryear]{elsarticle}
\usepackage[T1]{fontenc}
\usepackage[latin9]{inputenc}
\usepackage{geometry}
\usepackage{color}
\usepackage{array}
\usepackage{amsmath}
\usepackage{graphicx}
\usepackage{setspace}

\makeatletter

\providecommand{\tabularnewline}{\\}

\journal{xxx}

\makeatother

\usepackage{babel}
\begin{document}

\begin{frontmatter}

\title{Fuel treatment planning: fragmenting high fuel load areas while maintaining availability and connectivity of faunal habitat}

\author[rvt,focal]{Ramya Rachmawati}

\ead{ramya\_ unib@yahoo.co.id}

\author[rvt]{{\normalsize{}Melih Ozlen}  }

\ead{melih.ozlen@rmit.edu.au}

\author[rvt]{{\normalsize{}John Hearne} \corref{cor1} }

\ead{john.hearne@rmit.edu.au}

\author[rvt]{{\normalsize{}Karin Reinke}}

\ead{karin.reinke@rmit.edu.au}

\cortext[cor1]{Corresponding author}

\address[rvt]{School of Science, RMIT University, Melbourne, Australia}

\address[focal]{Mathematics Department Faculty of Mathematics and Natural Sciences,
University of Bengkulu, Bengkulu, Indonesia}

\begin{abstract}
\begin{doublespace}
{ Reducing the fuel load in fire-prone landscapes is aimed at mitigating the risk of catastrophic wildfires but there are ecological consequences. Maintaining habitat for fauna of both sufficient extent and connectivity while fragmenting areas of high fuel loads presents land managers with seemingly contrasting objectives. Faced with this dichotomy, we propose a Mixed Integer Programming (MIP) model that can optimally schedule fuel treatments to reduce fuel hazards by fragmenting high fuel  load regions while considering critical ecological requirements over time and space. The model takes into account both the frequency of fire that vegetation can tolerate and the frequency of fire necessary for fire-dependent species.  Our approach also ensures that suitable alternate habitat is available and accessible to fauna affected by a treated area. More importantly, to conserve fauna the model sets a minimum acceptable target for the connectivity of habitat at any time . These factors are all included in the formulation of a model that yields a multi-period spatially-explicit schedule for treatment planning. Our approach is then demonstrated in a series of computational experiments with  hypothetical landscapes,  a single vegetation type and a group of faunal species with the same habitat requirements. Our experiments show that it is possible to reduce the risk of wildfires while ensuring sufficient connectivity of habitat over both space and time. Furthermore, it is demonstrated that the habitat connectivity constraint is more effective than neighbourhood habitat constraints. This is critical for the conservation of fauna and of special concern for vulnerable or endangered species. }
\end{doublespace}
\end{abstract}

\begin{keyword}
Environmental modelling \sep Wildfire \sep Landscape management \sep Habitat conservation \sep  Optimisation 
\end{keyword}

\end{frontmatter}


\begin{doublespace}

\section{Introduction}
\end{doublespace}

\begin{doublespace}
Fire plays an important role in maintaining ecological integrity in many natural ecosystems \citep{Keane2010Evaluating} but wildfires also pose a risk to human life and economic assets \citep{King2008421}.  Climate change is expected to aggravate these risks \citep{kates2012transformational} but they can be reduced through  fuel management \citep{Agee200583,martell2015review, Ascoli2012920}. This is the process of altering the
structure and amount of  fuel accumulation in a landscape. It is achieved through the application
of treatments, such as prescribed burning or mechanical clearing.  To reduce the risk of large wildfires, fire management agencies in Australia \citep{McCaw2013217,Boer2009132}
and the USA \citep{Ager2010_FEM,Collins201024} have initiated extensive
fuel management programs in fire-prone areas.   Fuel load or biomass accumulation is a continuous ecosystem process. Each year parts of the landscape are treated to reduce the overall
fuel load for subsequent fire seasons.  Treatment frequency is partially dictated by the vegetation community.  Reducing the fuel load in the landscape in this way  helps to prevent or minimise the spread and intensity of wildfire.

Similarities exist between the fuel treatment problem described here
and the  planning problem for forest harvesting.
Both of these problems consider vegetation dynamics and can be seen
as a \textquoteleft timing problem\textquoteright , meaning that the
risk and values change over time as the vegetation grows. In the fuel
treatment problem, an area is treated to reduce fuel load; in the
forest harvesting problem, an area is harvested using mechanical clearing
for timber production.  Both activities have consequences for the habitat. 
Previous studies in the forest harvesting problem have taken into account some ecological
requirements. For example,  \cite{Bettinger1997111} used a Tabu search algorithm to schedule timber harvest subject to  spatial wildlife goals. Specifically, they maintained sufficient habitat of a certain maturity within a specified distance of a hiding or themal place. \citet{ohman2008incorporating} proposed an exact method
for long-term forest planning to maintain the biodiversity of the
forest. They believe that biodiversity in the forest ecosystem can
be maintained by minimising the total perimeter of old forest patches
so that the fragmentation of old forest is reduced. Hence,  compactness
of the habitat for species can be achieved. The model was run in a
five-yearly planning horizon across a landscape that comprised 924
stands . However, their model
did not consider habitat connectivity across time. Addressing this
shortcoming, \citet{konnyHu2014temporal} proposed a model that ensures
mature forest patches are temporarily connected between time-steps
while scheduling forest harvesting. The model achieves this without substantial reduction
in timber revenues. However,  this model does not take into account the overall
habitat connectivity of each period, nor does it track the habitat
connectivity across the entire planning horizon, both of which are
important for the persistence of species.

Various methods  have been proposed for incorporating the effect of wildfires into  harvest planning models. A comprehensive review is provided by \citet{Bettinger201043}.  More recently \citet{Troncoso2016} showed that including wildfire risks into a harvesting planning model with adjacency constraints can yield improved outcomes. The spatial arrangement of fuel treatment planning plays a substantial role in providing better
protection in the landscape \citep{Rytwinski2010}. Fuel arrangement
can modify fire behaviour and when fragmented, can lessen the chance
of large wildfires \citep{Kim2009253}. Considering the 'value at risk' \citet{Chung2013IJWF} used simulated annealing to  determine a long-term schedule for the location and timing of prescribed burns on a landscape.
An important factor that affects
wildfire extent is the connectivity of \textquoteleft old\textquoteright{}
untreated patches \citep{Boer2009132}. \citet{wei2014schedule} proposed
a single-period model to break the connectivity of high fuel load
patches by considering the duration and speed of a future
fire. Taking into account the vegetation dynamics over time is fundamental
to accurate fuel treatment planning \citep{Krivtsov20092915}. A multi-period model for fuel treatment planning that included the  dynamics of a single vegetation type was formulated by \citet{Minas2014412}. The objective in this model was to break the connectivity of \textquoteleft old\textquoteright{} patches
in the landscape  over the entire solution period of a few decades.

The efficacy of the applications of fuel treatment remains debated
among experts according to different perspectives \citep{penman2011prescribed}.
Fuel treatments reduce the overall fuel load in   landscapes 
but at the same time may result in significant habitat modification
for fauna living within the treated area. If the right mix of habitat
availability  in the landscape is not maintained, populations may be
adversely affected, leading to local extinctions where minimum viable
population thresholds are no longer met. For example, the Mallee emu-wren,
a native bird of Australia, depends on 15-year-old mallee-Triodia
vegetation for survival \citep{Brown2009432} . This vegetation recovers
very slowly after fuel treatments, and the Mallee emu-wren is unable
to survive in vegetation aged less than 15 years. Another Australian example is the Southern Brown Bandicoot. They require 5-15 year old heathland \citep{COBI:COBI934}.
Similarly, in California, frequent
fires can destroy the mature coastal sage scrub habitat
required for the coastal cactus wren and the California gnatcatcher
on which these species rely \citep{Conlisk201527predicting}. If we
want to conserve these species, it is important to maintain the availability
and connectivity of their habitats. In fact, more generally,    habitat connectivity is vital to support
the ecology and genetics of local populations \citep{rayfield2015multipurpose}.    
 The question then arises: Can fuel treatments be scheduled to break the connectivity of
high fuel load areas while maintaining the availability and connectivity of habitats?  

 Here we significantly extend current models by tracking and maintaining
defined levels of habitat connectivity over time, in addition to reducing
and fragmenting high fuel loads across the landscape. The model we
present is the first multi-period fuel treatment model that takes
into account habitat connectivity and solved using
exact optimisation. The proposed model is designed for  fire-dependent 
landscapes so additional ecological constraints are imposed based on the concept of Tolerable Fire Intervals (TFI's) \citep{Cheal2010}.  It is harmful for vegetation in an area to be subjected to another fire before a certain time (the minimum TFI) has elapsed   since the last fire in that area.  It is also desirable that a burn \textit{does} take place before a certain time ( the maximum TFI) has elapsed since the last fire.  Thus fuel treatment in each area is constrained to occur in a time-window between the minimum and maximum TFI since the last burn in that area. The TFI's are vegetation-dependent.

A Mixed Integer Programming (MIP) model is presented here for fuel
treatment planning. Subject to the time-windows imposed by the TFI's, the objective is 
to fragment  high fuel load areas as much as
possible while maintaining habitat connectivity in the landscape. The model is illustrated
with a single vegetation type and a single animal species. We assume that
the animals can relocate to a neighbouring area that has similar
habitat characteristics.
%
 The model is demonstrated on a series of hypothetical
landscapes.
\end{doublespace}

\begin{doublespace}

\section{Model formulation\label{sec:Model-formulation}}

\begin{figure}[h!]	
	\protect\caption{The relationship between the minimum TFI, mature, high fuel load,
		and the maximum TFI threshold values\label{fig:The-relationship-of_thresholds}}
	\begin{centering}
		\includegraphics[scale=0.35]{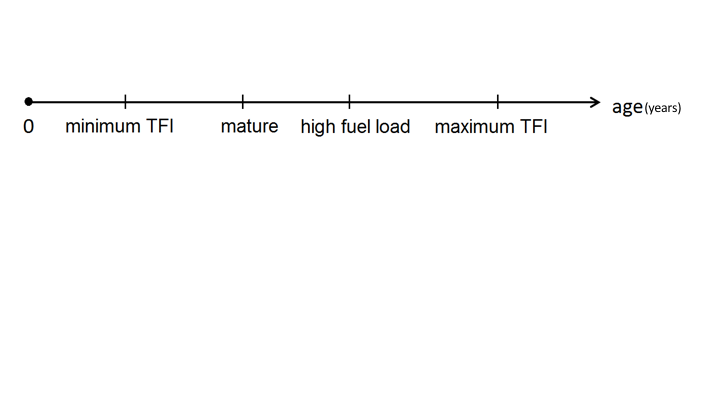}
		\par\end{centering}	
\end{figure}

In this formulation, cells represent the candidate locations for fuel treatment
in a landscape. The 'fuel age' (years) in each cell is defined as the time elapsed since the last treatment of that cell.
The cell\textquoteright s fuel age is reset to
zero if the cell is treated or incremented by one if untreated in any year. Each
cell has its minimum and maximum tolerable fire intervals (TFIs) which depend on the vegetation type in that cell. Within the time-window defined by the minimum and maximum TFI,  there is a time at which the vegetation is regarded as high risk from then on until the cell is treated. This time will be referred to as the 'high fuel load' threshold.  Without being specific, for this formulation we shall  consider a vertebrate that requires habitat offering 'mature' vegetation.  The vegetation age (time since last burnt) at which vegetation is considered 'mature'  will be referred to as the 'mature'  threshold.  In our example the mature threshold is less than the high fuel load threshold but the formulation is more general. The relationship between these thresholds is represented
in Figure \ref{fig:The-relationship-of_thresholds}. 

Reducing the connectedness of high fuel load cells through fuel treatment should reduce the risk  of fire spreading over a large area. Fuel treatment, however, modifies habitat.  For each cell treated in a given year, suitable neighbouring habitat should be available in the following year. Moreover,  for metapopulation persistence, to the extent possible we require the neighbouring habitat to be connected to other cells of mature habitat. 

\end{doublespace}

\selectlanguage{british}%
\begin{doublespace}
The following mixed integer programming model is formulated to determine a multi-period optimal schedule for treatment of 
cells.  The objective is to break the connectivity
of high fuel load cells in the landscape each year while providing continuity of habitat for the species of concern.

\textcolor{black}{\smallskip{}
}

\textcolor{black}{\bf Sets:}

\textcolor{black}{\smallskip{}
}

\textcolor{black}{$C$ is the set of all cells in the landscape}

\textcolor{black}{$\varPhi_{i}$ is the set of cells connected to
cell $i$}

\textcolor{black}{$T$ is the planning horizon}

\textcolor{black}{\smallskip{}
}

\textcolor{black}{\bf Indices:}

\textcolor{black}{\smallskip{}
}

\textcolor{black}{$i$ = cell}

\textcolor{black}{$t$ = year, $t$ = 0, 1, 2, \ldots$T$}

\textcolor{black}{\smallskip{}
}

\textcolor{black}{\bf Parameters:}

\textcolor{black}{\smallskip{}
}

\textcolor{black}{$a_{i}$ = initial fuel age of cell $i$}

\textcolor{black}{$R$ = the total area of }\foreignlanguage{british}{\textcolor{black}{cells
		in the landscape}}

$\rho$ = treatment level (percentage), i.e. the maximum fraction of $R$ that
can selected for treatment in any one year

$c_{i}$ = area of cell $i$

\selectlanguage{british}%
\textcolor{black}{$d_{i}$ = high fuel load  threshold for cell
}\textcolor{black}{\emph{i}}

\textcolor{black}{$m_{i}$ = mature  threshold for cell }\textcolor{black}{\emph{i}}

\textcolor{black}{$G_{t}$ = desired target of mature cell connectivity
in year }\textcolor{black}{\emph{t}}

\textcolor{black}{$MaxTFI_{i}$ = maximum tolerable fire interval
(TFI) of cell }\foreignlanguage{british}{\textcolor{black}{$i$}}

\textcolor{black}{$MinTFI_{i}$ = minimum TFI of cell }\foreignlanguage{british}{\textcolor{black}{$i$}}
	
\textcolor{black}{$M  \quad$ is a   "big M" parameter (must be greater than the maximum fuel age)}

\textcolor{black}{\smallskip{} }

\textcolor{black}{\bf Decision variables:}

\selectlanguage{british}%
\textcolor{black}{$A_{i,t}$ = fuel age of cell }\textcolor{black}{\emph{i}}\textcolor{black}{{}
in year }\textcolor{black}{\emph{t}}\textcolor{black}{{} }\foreignlanguage{english}{\textcolor{black}{\smallskip{}
}}

$x_{i,t}=\begin{cases}
1 & \mbox{if cell \emph{i} is treated in year \emph{t}}\\
0 & \mbox{otherwise}
\end{cases}$

\textcolor{black}{\smallskip{}
}

$Mature_{i,t}=\begin{cases}
1 & \mbox{if cell \emph{i} is classified as \textquoteleft mature \textquoteright\ in year \emph{t}}\\
0 & \mbox{otherwise}
\end{cases}$

\textcolor{black}{\smallskip{}
}

$HabitatConn_{i,j,t}=\begin{cases}
1 & \mbox{if connected cells \emph{i} and \emph{j} are both `mature' cells in year \emph{t}}\\
0 & \mbox{otherwise}
\end{cases}$

\textcolor{black}{\smallskip{}
}

$High{}_{i,t}=\begin{cases}
1 & \mbox{if cell \emph{i} is classified as high fuel load cell in year \emph{t}}\\
0 & \mbox{otherwise}
\end{cases}$

\textcolor{black}{\smallskip{}
}

$HighConn_{i,j,t}=\begin{cases}
1 & \mbox{if connected cells \emph{i} and \emph{j} are both high fuel load cells in year \emph{t}}\\
0 & \mbox{otherwise}
\end{cases}$

\textcolor{black}{\smallskip{}
}

$Old_{i,t}=\begin{cases}
1 & \textrm{if the fuel age in cell \emph{i }is classified as  over the maximum TFI in year } t \\
0 & \textrm{otherwise}
\end{cases}
$\textcolor{black}{\bigskip{}
}

\bigskip

{\bf The model}\\

\textcolor{black}{The objective  is to minimise $z$, the connectivity of high fuel load cells }
\end{doublespace}

\begin{equation}
\min z=\sum\limits^{T}_{t=1}\;\sum\limits_{i\in C} \;\sum\limits_{j\in\varPhi_{i},i<j}HighConn_{i,j,t} \label{eq:1}
\end{equation}

\selectlanguage{british}%
\begin{doublespace}
\textcolor{black}{subject to}
\end{doublespace}

\begin{singlespace}
\begin{equation}
\underset{i \in C}{\sum}c_{i}x_{i,t}\leq\rho R,\;t=1\ldots T    \label{eq:2}
\end{equation}
\begin{equation}
A_{i,0}=a_{i},\;\forall i\in C\label{eq:3}
\end{equation}

\begin{equation}
A_{i,t}\geq A_{i,t-1}+1-Mx_{i,t},\;t=1\ldots T,\forall i\in C\label{eq:4}
\end{equation}

\begin{equation}
A_{i,t}\leq M(1-x_{i,t}),\;t=1\ldots T,\forall i\in C\label{eq:5}
\end{equation}

\begin{equation}
A_{i,t}\leq A_{i,t-1}+1,\;t=1\ldots T,\forall i\in C\label{eq:6}
\end{equation}

\begin{equation}
A_{i,t}-d_{i}\leq M \: High_{i,t}-1,\;t=1\ldots T,\forall i\in C\label{eq:tr unit i is classified danger if age is d_=00007Bi=00007D or more_(1)}
\end{equation}

\begin{equation}
A_{i,t}\geq d_{i} \: High{}_{i,t},\;t=1\ldots T,\forall i\in C\label{eq:tr unit i is classified danger if age is d_=00007Bi=00007D or more_(2)}
\end{equation}

\begin{equation}
High_{i,t}+High{}_{j,t}- HighConn_{i,j,t}  \leq1,\;t=1\ldots T,\forall j\in\varPhi_{i},i<j,\forall i\in C\label{eq:connectivity of treatment unit}
\end{equation}

\begin{equation}
A_{i,t}-m_{i}\leq M \: Mature_{i,t}-1,\;t=1\ldots T,\forall i\in C\label{eq:tr unit i is classified mature if age is m_=00007Bi=00007D or more_(1)}
\end{equation}

\begin{equation}
A_{i,t}\geq m_{i} \: Mature_{i,t},\;t=1\ldots T,\forall i\in C\label{eq:tr unit i is classified mature if age is m_=00007Bi=00007D or more_(2)}
\end{equation}

\begin{equation}
\underset{j\in\varPhi_{i}}{\sum}Mature{}_{j,t}\geq x_{i,t},\;t=1\ldots T,\forall i\in C\label{eq:if there is at least one neighbouring mature cell next year, then we can burn cell}
\end{equation}

\begin{equation}
Mature{}_{i,t}+Mature{}_{j,t}-HabitatConn_{i,j,t}\leq1,\;t=1\ldots T,\forall j\in\varPhi_{i},i<j,\forall i\in C\label{eq:HabitatConn1}
\end{equation}

\begin{equation}
Mature{}_{i,t}+Mature{}_{j,t}\geq2HabitatConn_{i,j,t},\;t=1\ldots T,\forall j\in\varPhi_{i},i<j,\forall i\in C\label{eq:HabitatConn1-1}
\end{equation}

\begin{equation}
\underset{i\in C}{\sum} \sum\limits_{j\in\varPhi_{i}, i<j}HabitatConn{}_{i,j,t}\geq G_{t},\:t=1\ldots T   \label{eq:28-1}
\end{equation}

\begin{equation}
A_{i,t}-MaxTFI_{i}\leq M \: Old_{i,t}-1,\;t=0\ldots T-1,\forall i\in C\label{eq:tr unit i is classified 'over max TFI' (1)}
\end{equation}

\begin{equation}
A_{i,t}\geq MaxTFI_{i} \: Old{}_{i,t},\;t=0\ldots T-1,\forall i\in C\label{eq:tr unit i is classified 'over max TFI' (2)}
\end{equation}

\begin{equation}
Old{}_{i,t-1}+\frac{1}{\mid\varPhi_{i}\mid}\underset{j\in\varPhi_{i}}{\sum}Mature{}_{j,t}\leq1+x_{i,t},\;t=1\ldots T,\forall i\in C\label{eq:if cell i is over max TFI burn if there is at least one neighbouring mature cell next year. But don't burn if there's no neighbouring mature cell-1}
\end{equation}

\begin{equation}
A_{i,t-1}\geq MinTFI_{i} \: x_{i,t},\;t=1\ldots T,\forall i\in C\label{eq:dont_burn_if_below_minTFI}
\end{equation}

\begin{equation}
x_{i,t},\;High{}_{i,t},\;HighConn_{i,j,t},\;Mature{}_{i,t},\;Old{}_{i,t}\in\{0,1\}\label{eq:binary vars}
\end{equation}

\end{singlespace}

\begin{doublespace}
\textcolor{black}{The objective function \eqref{eq:1} minimises the
connectivity of high fuel load cells in a landscape across the planning
horizon. Constraint \eqref{eq:2}} specifies
that the total area selected for fuel treatment each year should not exceed
a fixed proportion of the totalarea of the landscape.
\textcolor{black}{Constraint \eqref{eq:3} sets the initial fuel age
in a cell. Constraints \eqref{eq:4} to \eqref{eq:6} track the fuel
age of each cell. 
 Constraints \eqref{eq:4} and \eqref{eq:6}
increment fuel age by exactly one year if the cell is not treated. Constraint \eqref{eq:5} forces the fuel age to be reset to zero if the cell is treated. Note that the $A_{i,t}$ are continuous variables although only integer values are assigned to them. }

Constraints \eqref{eq:tr unit i is classified danger if age is d_=00007Bi=00007D or more_(1)}
and \eqref{eq:tr unit i is classified danger if age is d_=00007Bi=00007D or more_(2)}
use binary variable \textcolor{black}{$High{}_{i,t}$ to classify
a cell as a high fuel load cell if and only if the fuel age exceeds
a threshold value}. In Constraint \eqref{eq:connectivity of treatment unit},
$HighConn_{i,j,t}$ takes the value one if connected cells \emph{i}
and \emph{j} are both classified as high fuel load cells in year\emph{
t}.

\selectlanguage{british}%
\textcolor{black}{Constraints \eqref{eq:tr unit i is classified mature if age is m_=00007Bi=00007D or more_(1)}
to \eqref{eq:tr unit i is classified mature if age is m_=00007Bi=00007D or more_(2)}
classify a cell to be a } `mature'
cell, if and only if the fuel age is over the mature age threshold.
Constraint \eqref{eq:if there is at least one neighbouring mature cell next year, then we can burn cell}
states that we cannot treat a cell in this period unless there is at
least one neighbouring mature cell in the following year.

In this model, we also ensure that sufficient habitat (mature-cell)
connectivity in the landscape as a whole is available each year. Constraints \eqref{eq:HabitatConn1}
and  \eqref{eq:HabitatConn1-1} ensure that $HabitatConn{}_{i,j,t}$
takes the value one if and only if connected cells $i$ and $j$
are both classified as mature cells in year $t$. Constraint \eqref{eq:28-1}
ensures that the number of habitat connections each year is greater
than the desired target, $G_{t}$. 

\textcolor{black}Constraints \eqref{eq:tr unit i is classified 'over max TFI' (1)}
to \eqref{eq:tr unit i is classified 'over max TFI' (2)} classify
a cell as 'Old'  if and only if the fuel age is over the maximum TFI. Constraint
\eqref{eq:if cell i is over max TFI burn if there is at least one neighbouring mature cell next year. But don't burn if there's no neighbouring mature cell-1}
ensures that a cell must be treated if the cell\textquoteright s fuel
age is over maximum TFI, and there is at least one neighbouring mature
cell in the following period. This constraint avoids a deadlock that
may occur when the cell's fuel age is over the maximum TFI and there
are no neighbouring mature cells for the next period. In this study,
we break the deadlock in favour of mature cell availability. Constraint
\eqref{eq:dont_burn_if_below_minTFI} ensures that the cell with fuel
age less than the minimum TFI cannot be treated. Constraint
\eqref{eq:binary vars} ensures that the decision variables take binary
values.

\end{doublespace}

\begin{doublespace}

\section{Model illustration}
\end{doublespace}

\begin{figure}
\protect\caption{Initial proportion of cells in the landscape of each fuel age group
for the computational
experiments\label{fig:Proportion-of-initial}}

\medskip{}

\centering{}\includegraphics[scale=0.5]{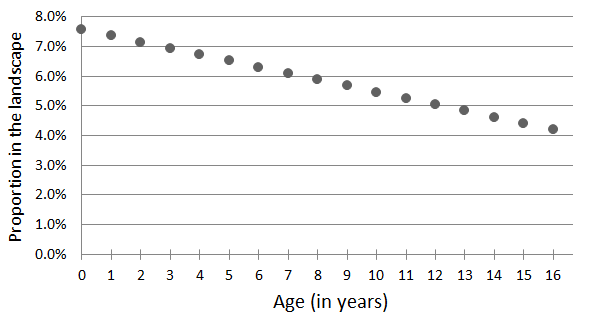}
\end{figure}

\begin{figure}[h]
\protect\caption{The definition of connected cells. Cell 5 is considered connected
to cells 6 (right) , 4 (left), 2 (up) and 8 (down)\label{fig:The-connected-cells}}

\begin{centering}
\includegraphics[scale=0.5]{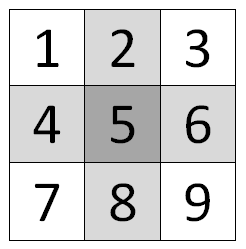}
\par\end{centering}

\end{figure}

\begin{doublespace}
In this section, we demonstrate the model formulated in Section
\ref{sec:Model-formulation} using hypothetical random landscapes
comprising 100 grid cells, generated using the NLMpy package \citep{Etherington2015NLMpy}.
(Note that the model does not require a regular grid. Cells can be any shape and all that is needed is that the neighbours of each cell are specified.) For this illustration we assume that there is a single fuel type in the landscape, with
the thresholds of mature (suitable habitat) and  high fuel load ages set as 8 and
12 years old, respectively. The minimum and the maximum TFIs are chosen
as 2 and 16 years. The initial fuel ages in the landscape
are between 0 and 16 years, this means that not all the cells are
categorised as high fuel load. Figure \ref{fig:Proportion-of-initial}
represents the assumed distribution of the initial cell fuel age.
For this illustration a cell is assumed to be connected to its immediate neighbouring cells
that have shared boundaries (Figure \ref{fig:The-connected-cells}).
Suppose that there are at most ten cells to be treated each year (ten
percent of the total area in the landscape), and the length of planning
horizon is 13 years. 
\end{doublespace}

\begin{figure}
\protect\caption{Illustration of initial high fuel load cell and habitat connectivity
in the landscape, the arrow (\ensuremath{\leftrightarrow}) represents
one connection\label{fig:Initial-high-risk-cell}}

\centering{}\includegraphics[scale=0.45]{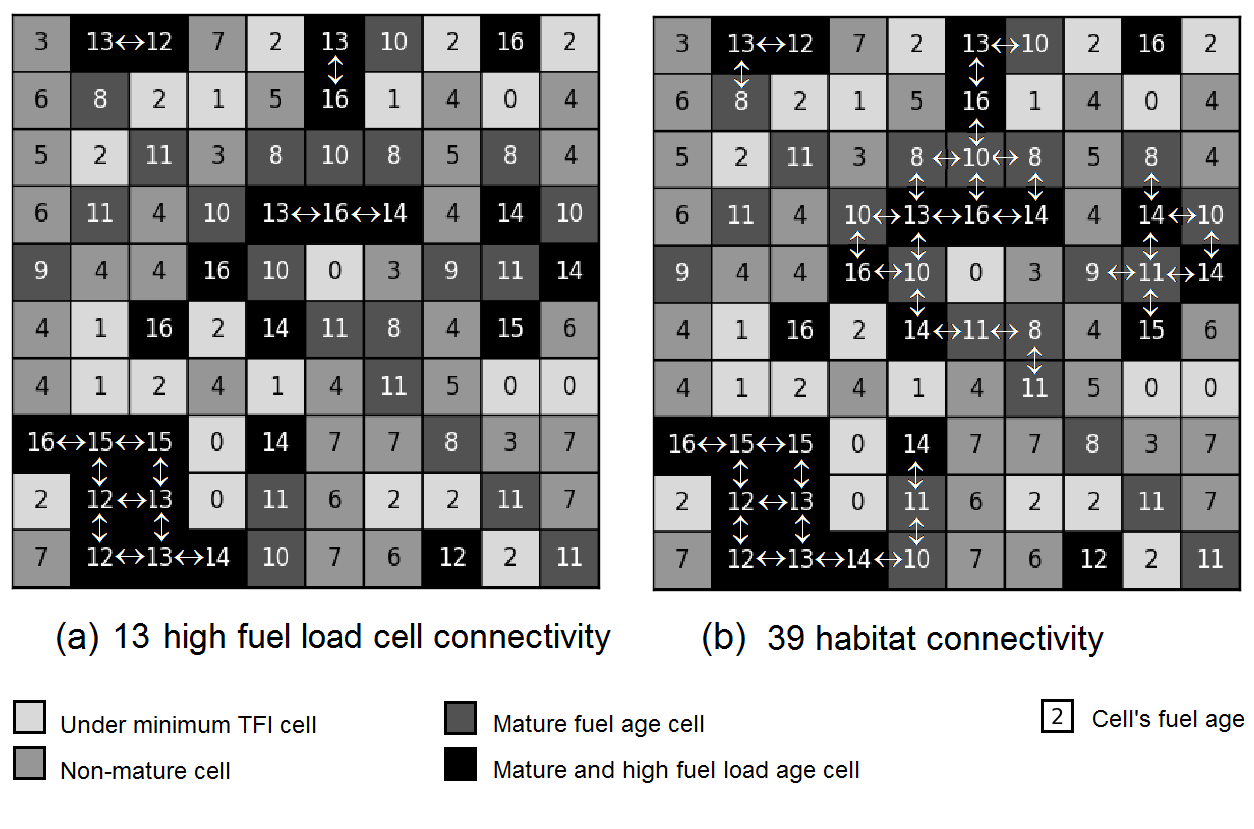}
\end{figure}

\begin{doublespace}
As shown in Figure \ref{fig:Initial-high-risk-cell},  initially the landscape has 13 high fuel load cell connections that we want to reduce with time. It also has 
39 habitat connections that we want to maintain over the planning
horizon.
In this model illustration, we compare four different settings (Table
\ref{tab:Four-scenarios-Table}).

\begin{table}[h!]
	\protect\caption{Four settings for the model illustration and the computational experiments\label{tab:Four-scenarios-Table}}
		
	\medskip{}	
	
	\centering{}%
	\begin{tabular}{|c|c|>{\centering}m{3cm}|}
		\hline 
		& {\footnotesize{}Limit on overall habitat connectivity} & {\footnotesize{}Neighbouring habitat cell requirement for treatment}\tabularnewline
		\hline 
		{\footnotesize{}Setting 1} & {\footnotesize{}$G_{t}$ is set to the initial number of habitat connectivity
			of the landscape} & {\footnotesize{}yes}\tabularnewline
		\hline 
		{\footnotesize{}Setting 2} & {\footnotesize{}$G_{t}$ is set to the initial number of habitat connectivity
			of the landscape} & {\footnotesize{}no}\tabularnewline
		\hline 
		{\footnotesize{}Setting 3} & {\footnotesize{}$G_{t}$ is set to zero} & {\footnotesize{}yes}\tabularnewline
		\hline 
		{\footnotesize{}Setting 4} & {\footnotesize{}$G_{t}$ is set to zero} & {\footnotesize{}no}\tabularnewline
		\hline 
	\end{tabular}
\end{table}

 In the first and second settings,
we maintain the initial habitat connectivity, at a minimum
level of 39 connections. In the first setting we enforce the requirement
that a cell can only be treated if there is a neighbouring cell forming
a suitable habitat, but in the second setting that
requirement is relaxed. In the third setting,  only the neighbouring habitat
cell requirement is enforced without maintaining the overall habitat
connectivity. Setting 4 represents the base case with the only aim
of fragmenting high fuel load cells without habitat considerations. 

\end{doublespace}

\begin{doublespace}
\section{Illustration Results}

\begin{figure}
	\protect\caption{Fuel treatment schedule with ten percent treatment level and thirteen-year
		planning horizon for the first setting, $G_{t}$= initial\label{fig:Burn-schedule-Sc_1}}
	\medskip{}
	\centering{}\includegraphics[scale=0.4]{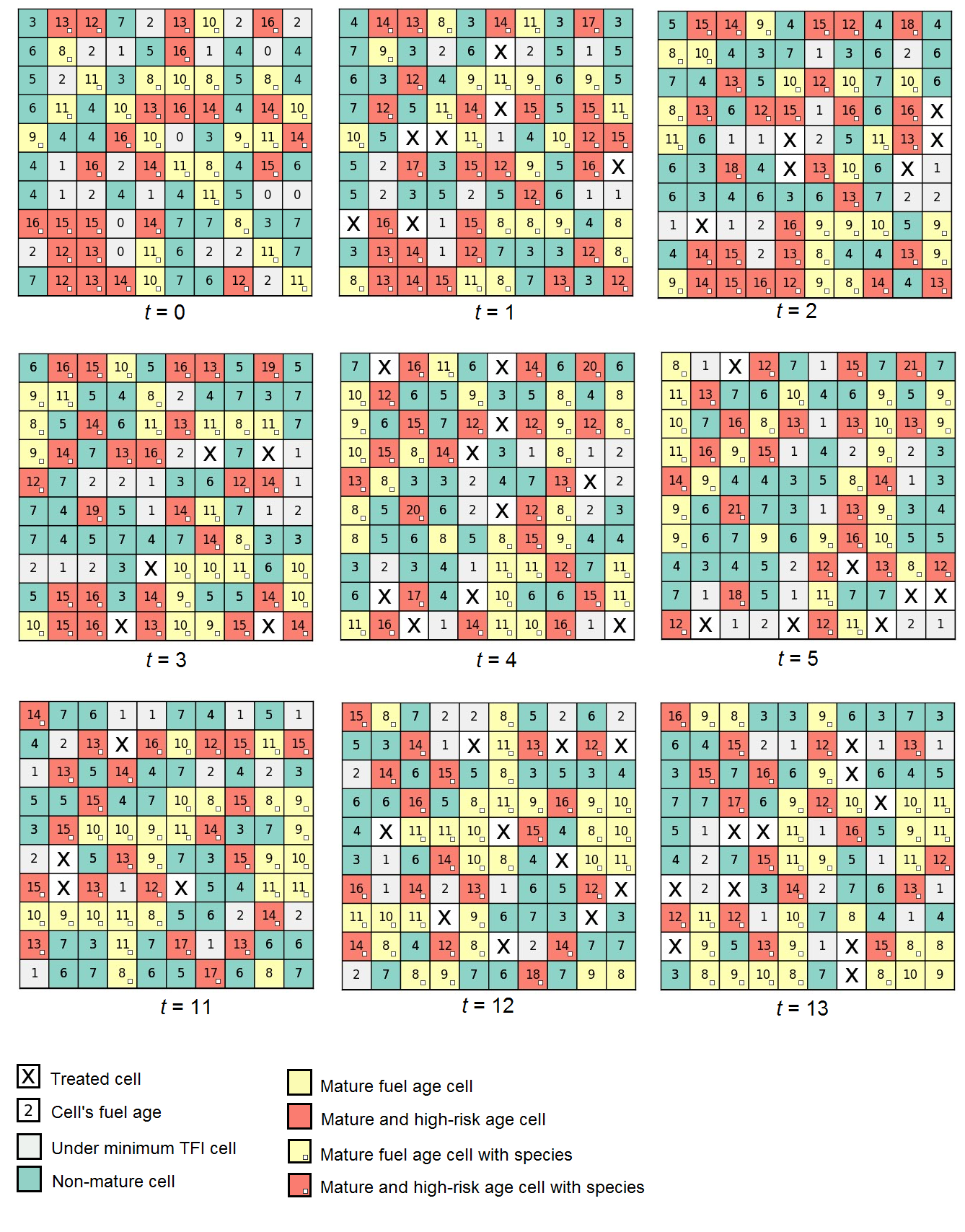}
\end{figure}

A sequence of landscape mosaics for the solution to setting 1 is given in Figure \ref{fig:Burn-schedule-Sc_1}. At $t=0$ note that the fuel age of cell $(1,9)$ has reached its maximum TFI. It is not selected for treatment as there is no neighbouring cell with suitable habitat i.e. no neighbouring mature cell (age $\ge$ 8). Recall that, for this illustration, only cells that share a common boundary are regarded as neighbours. Thus even at $t=5$ this cell is not considered for treatment. At this stage, however, the two row neighbours both have a fuel age of 7 and so at $t=6$ will provide suitable 'mature' habitat and the cell is in fact treated at this time (not shown but can be deduced from the fuel age shown at $t=11$).

It is also worth noting the four cells in the bottom right hand corner. Initially two of these cells are occupied. At $t=5$ the animals have moved to suitable neighbouring habitat and the cells remain unoccupied until $t=13$ when recolonisation has begun. It appears that the model is achieving the conservation goals. It is easy to see that the fragementation of high fuel cells has also been achieved. None of the high fuel load cells (in red) have a high fuel load neighbour.

\begin{figure}
\protect\caption{Habitat connectivity and high fuel load connectivity
 for the illustrative example} \label{fig:The-comparison-of}

\centering{}\includegraphics[scale=0.45]{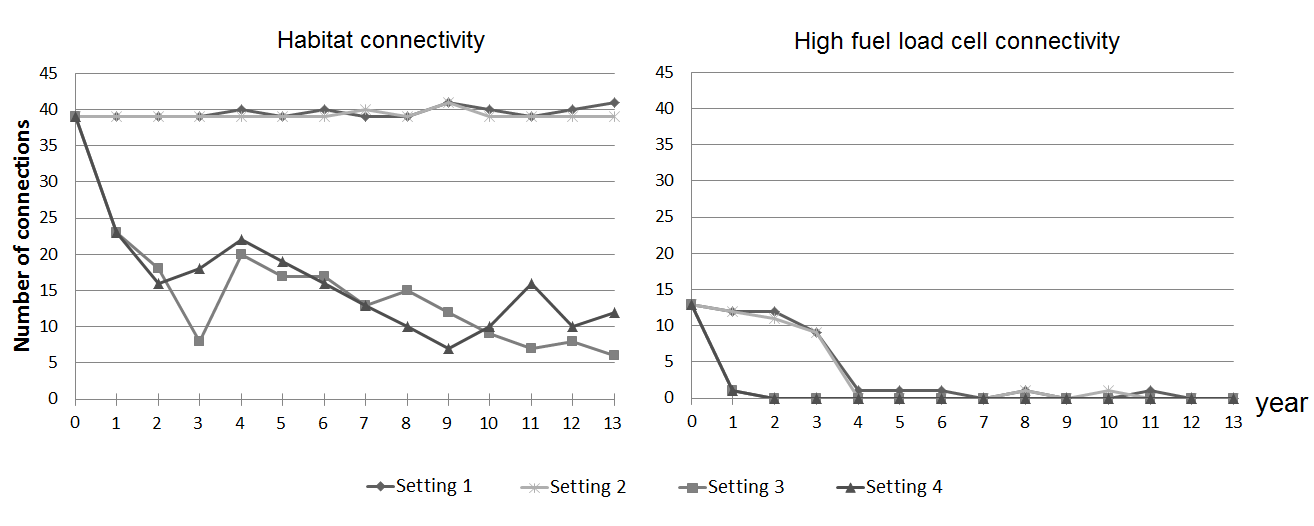}
\end{figure}

\begin{figure}
\protect\caption{The percentages of high fuel load cells and mature cells in the landscape
for the model illustration \label{fig:The-comparison-of-prop_model_demonstration}}

\centering{}\includegraphics[scale=0.4]{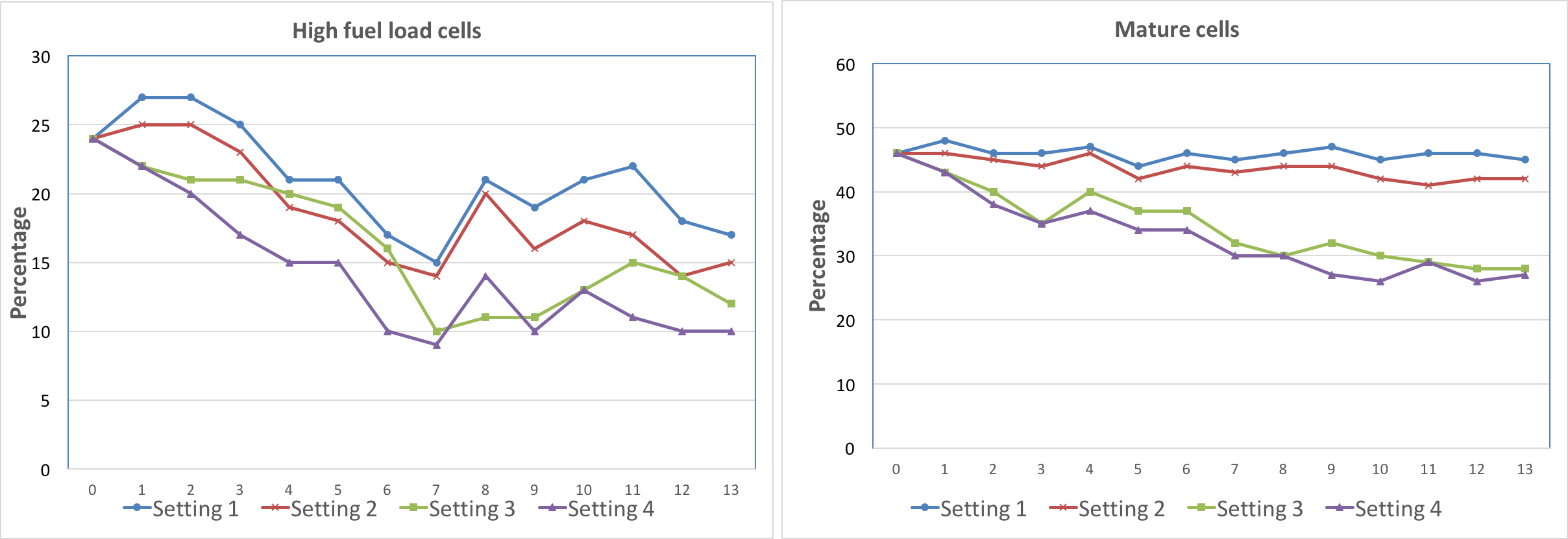}
\end{figure}

A comparison of the results  for all four settings is shown in  Figure \ref{fig:The-comparison-of}.  The high fuel load cells in the landscape are fully fragmented more quickly for settings 3 and 4 than settings 1 and 2. This is to be expected as the habitat constraints are relaxed for settings 3 and 4 and habitat connectivity drops rapidly as a consequence. Nevertheless, from $t=4$ on settings 1 and 2 do achieve similar fuel fragmentation while maintaining habitat connectivity throughout. In this case, however, Figure \ref{fig:The-comparison-of-prop_model_demonstration} shows that the landscapes comprise a greater number of high fuel load cells. On the other hand settings 3 and 4 not only perform poorly with regard to habitat \textit{connectivity} but habitat availability (mature cells) also declines as seen in Figure \ref{fig:The-comparison-of-prop_model_demonstration}.

\end{doublespace}

\begin{figure}
\protect\caption{The percentage of mature cells in the landscaape with animals present. 
\label{fig:Proportion-of-mature-and-species-model-illustration}}

\centering{}\includegraphics[scale=0.45]{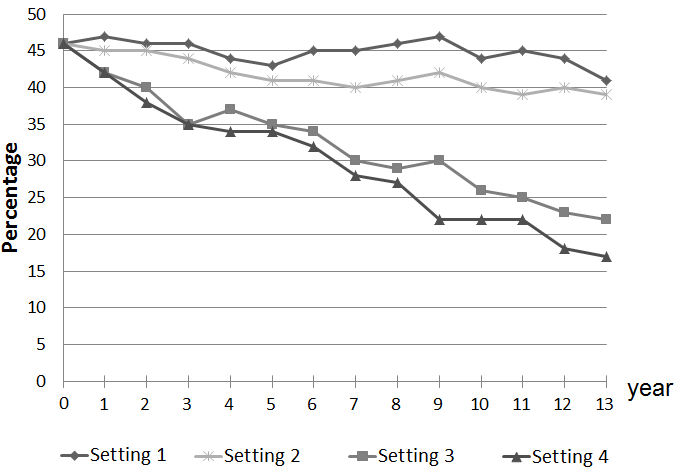}
\end{figure}

\begin{doublespace}
The location  of animals in mature cells in the landscape can be tracked over the planning horizon for the four settings.  It is assumed that initially all mature cells (includes high
fuel load cells) are populated by a particular vertebrate of interest. In any given year the vertebrate can only move to a neighbouring cell with suitable habitat. If a cell that is treated has no suitable neighbouring cell then any animals in that cell will die. An unoccupied mature cell can be (re)colonised from an occupied neighbour. An analysis for the four settings can be undertaken using  Figure
\ref{fig:Burn-schedule-Sc_1} for setting 1 and similar graphs (not shown) for the other three settings. The results are shown in Figure
\ref{fig:Proportion-of-mature-and-species-model-illustration}. The value of the connectivity constraints is now even more apparent. By the end of the planning horizon only 17\% of the landscape is occupied in setting 4 as opposed to 41\% for setting 1. Furthermore, setting 2 which includes the connectivity constraint but not the neighbourhood constraint ends up with 39\% of the landscape occupied compared with only 22\% for setting 3. Recall that setting 3  imposes neighbourhood but not connectivity constraints.


\end{doublespace}

\section{Computational experiments}

\begin{doublespace}
The illustration in the previous section was for a particular configuration of initial fuel age of cells in a 10x10 landscape. Were the previous findings simply a consequence of the initial configuration? In this section we explore landscapes with randomly generated initial configurations but with the same proportions of initial fuel age cells as given in  Figure \ref{fig:Proportion-of-initial}. 

We consider landscape sizes of 10x10 and 15x15 cells. In each case 30 landscapes were generated using the NLMpy package.  The model was solved for each of the four settings given in Table \ref{tab:Four-scenarios-Table}. A ten percent treatment level was applied with a planning horizon of 10 years. For the first two settings, we evaluated the initial number of connected mature cells for each landscape. This value of habitat connectivity, $G_t$, was then maintained over the planning horizon by constraint (15). We found, however, that for some landscapes it is impossible to maintain the initial extent of habitat over the
planning horizon. To deal with this infeasibility, we ran the model
by assigning a lower value of $G_{t}$ for the first years of a planning
horizon, and setting the higher value (the initial level of habitat
connectivity) of $G_{t}$ for the remainder of the planning
horizon only once it was feasible. 
 
The computational experiments were conducted using ILOG
CPLEX 12.6.2 with the Python 2.7.2 programming language using PuLP
modeller. The experiments were ran on Trifid, a computer cluster of
V3 Alliance. A single node with 16 cores of Intel Xeon E5-2670 and 64
GB of RAM was used.
 
\end{doublespace}

\begin{figure}
\protect\caption{High fuel load connectivity and habitat
connectivity with 95\% confidence intervals for the computational experiments \label{fig:95=000025-confidence-interval}}

\centering{}\includegraphics[scale=0.4]{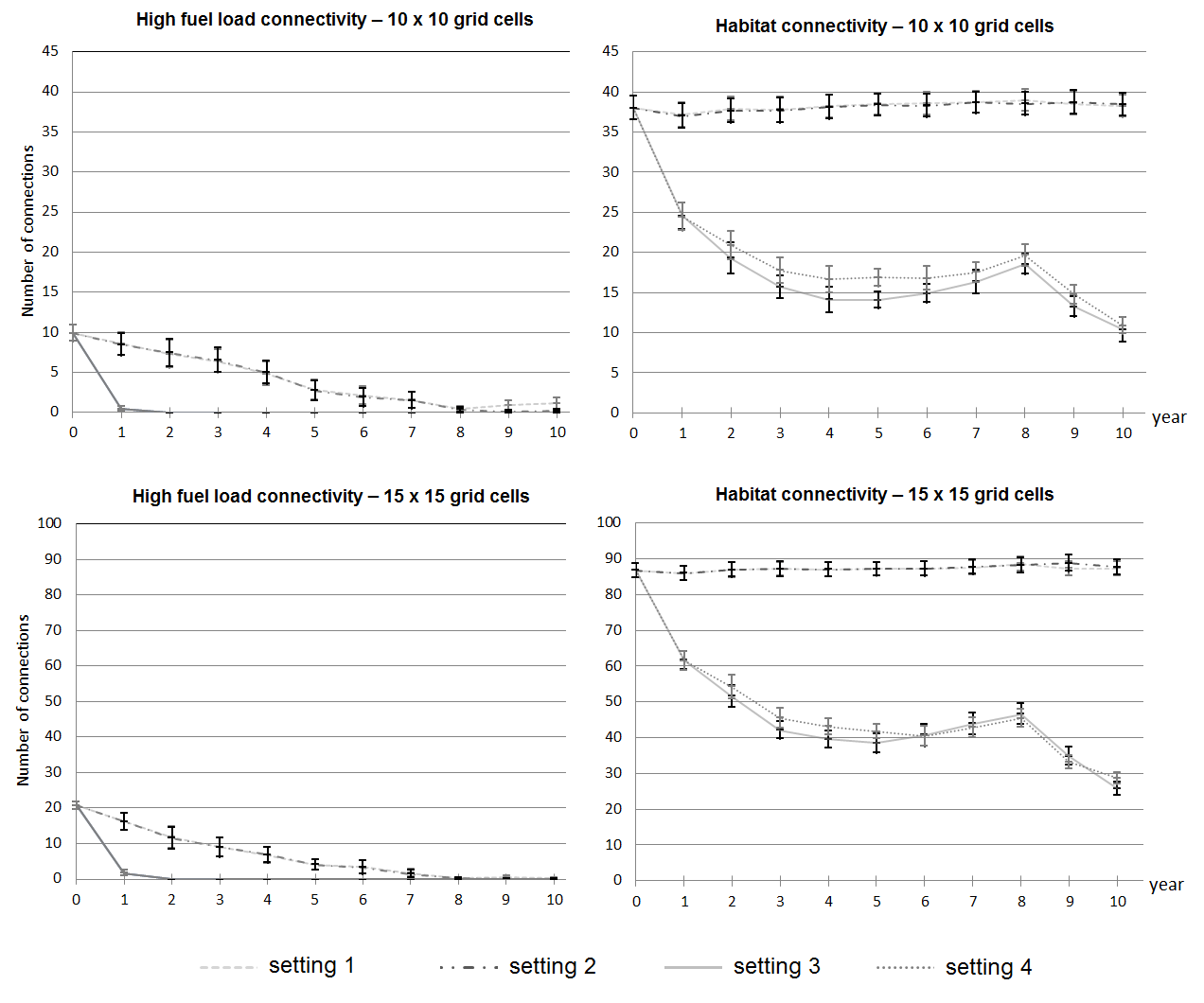}
\end{figure}

\begin{figure}
\protect\caption{Proportions of high fuel load cells and
mature cells in the landscape with 95\% confidence intervals for the computational experiments\label{fig:95=000025-confidence-interval-prop}}

\centering{}\includegraphics[scale=0.45]{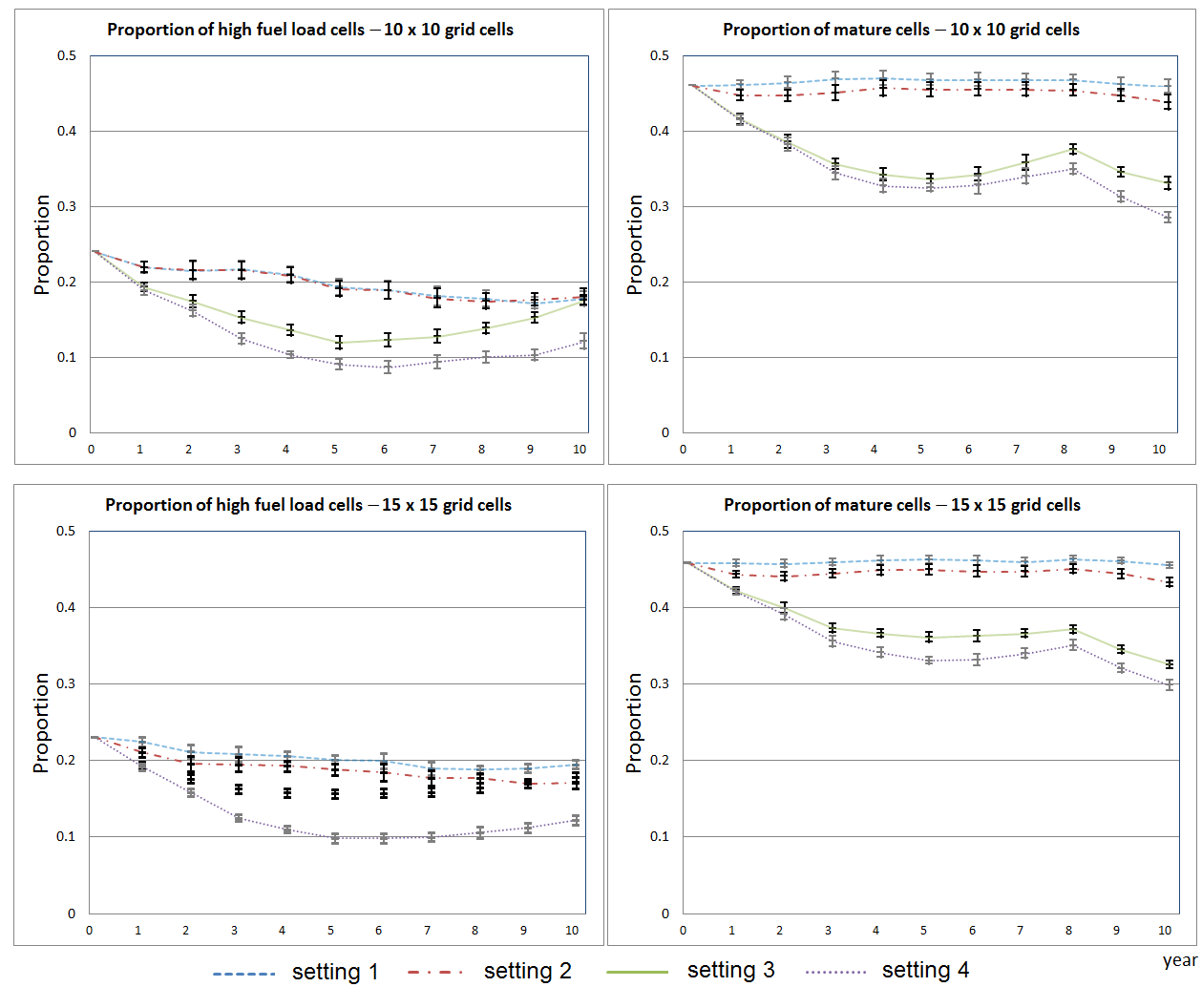}
\end{figure}

\begin{figure}
\protect\caption{Proportion of mature cells in the landscape with a faunal presence  for the computational
experiments\label{fig:Proportion-of-mature-comp-ex}}

\centering{}\includegraphics[scale=0.5]{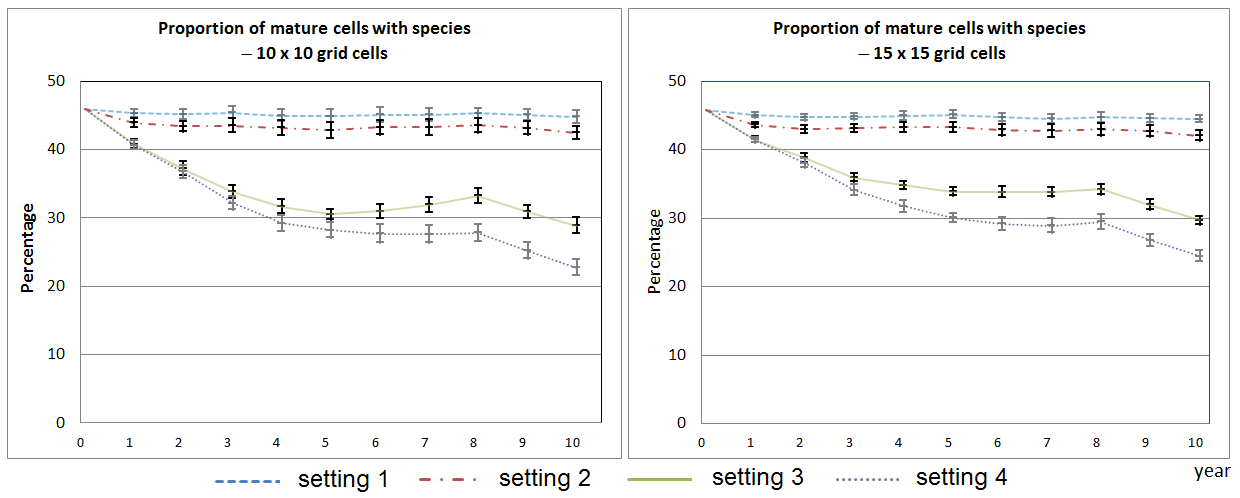}
\end{figure}

\section{Results of computational experiments}

\begin{doublespace}

Overall this more comprehensive analysis does not reveal any surprising differences from that observed in the model illustration. Figure \ref{fig:95=000025-confidence-interval} shows that, on average, settings 1 and 2  do reduce the high fuel load connectivity  but more slowly than in the case of the model illustration. On the other hand settings 3 and 4 achieve a rapid reduction in high fuel load connectivity but to the detriment of habitat connectivity. Figure \ref{fig:95=000025-confidence-interval-prop} shows, not unexpectedly, that
 settings 1, 2 and 3 all leave a greater proportion of high fuel load cells in the landscape compared with setting 4. Given that the difference betweeen setting 2 and setting 3 is that the former is concerned only with habitat connectivity and the latter with only with suitable neighbouring sites  Figure \ref{fig:Proportion-of-mature-comp-ex} reveals
  a remarkable difference between their performance in maintaining sites with a faunal presence.

\end{doublespace}

\begin{doublespace}
\section{Discussion}

The model results show that it is possible to achieve reductions in the number and  connectivity of high fuel load cells in the landscape while simultaneously ensuring habitat indices are maintained at their initial levels. While the reduction in the number of high fuel load cells is not as good when  habitat connectivity constraints are imposed, even in this case there is still a significant reduction in the overall {\it connectivity} of high fuel load cells. The work of \citet{wei2014schedule} indicates that this fragmentation of high fuel load areas is likely to reduce the risk of large  wildfires. 

Two methods were considered to meet conservation goals. One method was to maintain habitat connectivity and the other was to ensure that no cell was treated unless there was suitable habitat in its neighbourhood. This latter method is closely related to considering suitable habitat or forage in the neighbourhood of hiding places for a vertebrate (see for example \cite{Bettinger1997111}. The results clearly suggest, however, that the habitat connectivity contraints we used for setting 1 and 2 produced a significantly better outcome in terms of the fraction of the landscape occupied by our representative faunal species.

 In our model we defined the neighbourhood set to be the same for both high fuel load as well as habitat.  In practice, the set of habitat cells and the set of high fuel load cells forming the neighbourhood of a given cell will differ. In the case of a high fuel load cell the neighbouring cells could be weighted to take fire spread dynamics into account. In the case of habitat, neighbouring cells would need to be defined in terms of the  particular requirements and mobility of denizens living in a cell selected for treatment. In both cases of high fuel load and habitat, 'neighbourhoods' might comprise  more than just adjacent cells. Mathematically, this is easy to accommodate. The sets $\Phi_i$ used in constraints \eqref{eq:HabitatConn1}, \eqref{eq:HabitatConn1-1}  and \eqref{eq:if cell i is over max TFI burn if there is at least one neighbouring mature cell next year. But don't burn if there's no neighbouring mature cell-1} would simply be replaced by another set $\Psi_i$, say, specifying the sites that form the neighbourhood of site $i$.

The model presented in this paper comprises hypothetical landscapes with a single vegetation type  and a single faunal species. The model was developed  particularly for a fire-dependent vegetation type in a fire-prone landscape. An extension of the model to multiple vegetation types without the habitat connectivity has already been demonstrated on a real landscape \citep{Rachmawati201694}.  In principle, extensions to include multiple groups of faunal species can be achieved with the inclusion of additional constraints of the type  \eqref{eq:HabitatConn1}, \eqref{eq:HabitatConn1-1}  and \eqref{eq:if cell i is over max TFI burn if there is at least one neighbouring mature cell next year. But don't burn if there's no neighbouring mature cell-1}.  In practice habitat connectivity would need to be limited to a few groups of species. The needs of keystone species and vulnerable or endangered species would require particular attention. To some extent the problem is a  dynamic version of the Reserve Design Problem  \citep{Wang2016413,jafari2013new}. In this case the landscape from which areas for the reserve are to be selected change each year. Moreover decisions made in one period affect the subsequent landscape and hence the actions to be taken in future periods.

\end{doublespace}

\begin{doublespace}

\section{Conclusion}
\end{doublespace}

\begin{doublespace}

In this paper, we proposed and tested a mixed integer programming model that aimed to simultaneously fragment areas of high fuel load while maintaining the initial level of habitat connectivity.  The model was tested on a hypothetical landscape with a single vegetation type and a single faunal species with the same habitat needs. Some reduction in high fuel load areas could still be achieved after imposing a habitat connectivity constraint. Perhaps more importantly it was possible to achieve  significant overall reductions in high fuel load connectivity while maintaining habitat connectivity. The model, designed for fire-dependent landscapes achieves these outcomes whilst also ensuring that the vegetation is subject to fire of a necessary and sufficient frequency within tolerable limits.

The approach was  based on a theoretical perspective and has not
yet been applied to real landscapes. Nevertheless, a model based on a similar concept with mulitple vegetation types but without the habitat connectivity considerations has been successfully applied to a real landscape and closely related problems in harvest planning  have successfuly applied heuristics such as simulated annealing and Tabu search.

The development of optimised solutions for conflicting objectives has the potential to improve planning and operational decision making of prescribed burning strategies.  It is hoped that our approach can assist fire and land management agencies in making
their decisions about the timing and locations of future fuel treatments
while considering  critical ecological requirements.  For this purpose we
plan to extend the model to include multiple types of habitat
and species in the landscape.

Spatial optimisation models addressing 'connectivity' have been developed before for various purposes.  Such models are useful, for example,  when parcels of land need to be acquired for a particular purpose such as an airport or golf course. In the problem  addressed here not only is 'connectivity' in the landscape required for one attribute but also required in the same landscape is disconnectedness or fragmentation for another attribute. This new class of problem could prove useful for other purposes such as designing a connected reserve for an endangered species while  fragmenting the habitat needs for an invasive species.

\end{doublespace}

\begin{doublespace}

\section*{\textcolor{black}{Acknowledgment}}
\end{doublespace}

\begin{doublespace}
\textcolor{black}{The first author is supported by the The Directorate
General of Resources for Science, Technology and Higher Education,
the Ministry of Research, Technology and Higher Education of Indonesia
(1587/E4.4/K/2012). The second author is supported by the Australian
Research Council under the Discovery Projects funding scheme (project
DP140104246).}
\end{doublespace}

\begin{doublespace}

\section*{References}
\end{doublespace}

\begin{doublespace}
\bibliographystyle{elsarticle-harv}
\bibliography{References}
\end{doublespace}

\end{document}